\documentclass{article}
\usepackage{amsmath,amssymb,amsthm,amsfonts,latexsym,graphicx,mathtools}

\usepackage{url}
\usepackage{xcolor}

\DeclareMathOperator{\inp}{in}
\DeclareMathOperator{\out}{out}
\DeclareMathOperator{\id}{Id}

\usepackage{bbm}
\usepackage{enumitem}

\usepackage{comment}

\theoremstyle{plain}
\newtheorem{theorem}{Theorem}[section]
\newtheorem*{theorem*}{Theorem}
\newtheorem{lemma}[theorem]{Lemma}

\newtheorem{corollary}[theorem]{Corollary}

\theoremstyle{definition}

\theoremstyle{remark}
\newtheorem{remark}[theorem]{Remark}

\usepackage[affil-it]{authblk}
\usepackage{booktabs}   

\usepackage{microtype}

\title{The real eigenvalues of the non-backtracking Laplacian}
\author[1]{Raffaella Mulas\thanks{r.mulas@vu.nl}}
\author[2]{Leo Torres\thanks{leo@leotrs.com}}
\affil[1]{Vrije Universiteit Amsterdam, Amsterdam, The Netherlands}
\affil[2]{Independent researcher, Leipzig, Germany}

\date{}

\begin{document}

\maketitle

\begin{abstract}

Let $G$ be a simple graph on $N$ nodes and $M$ edges, with minimum degree at least $2$. We prove that the non-backtracking Laplacian of $G$ has at least $2M-2N+2c+2b$ real eigenvalues, counted with algebraic multiplicity, where $c$ is the number of connected components of $G$ and $b$ is the number of bipartite ones. More precisely, at least $M-N+2c$ of them lie in $[0,1)$ and at least $M-N+2b$ lie in $(1,2]$. We show that the bound $2M-2N+2c+2b$ is sharp, and that a connected $d$-regular graph attains it if and only if it is strictly Ramanujan, i.e., every adjacency eigenvalue $\mu\ne\pm d$ satisfies $|\mu|<2\sqrt{d-1}$. For completeness, we also show that the non-backtracking Laplacian has at most $2M-2c-2b$ real eigenvalues.

\vspace{0.2cm}
\noindent {\bf Keywords:} non-backtracking random walks; spectral graph theory; non-backtracking matrix; non-backtracking Laplacian; Ramanujan graphs
\end{abstract}

\section{Introduction}

While most of the spectral theory of simple graphs concerns symmetric matrices, such as the adjacency matrix and the graph Laplacians (we refer to \cite{chung,brouwer,BookJMZ} for three monographs), a growing body of work is devoted to non-symmetric operators, and in particular to non-backtracking operators. The eigenvalues of these operators are in general complex, and even a basic question such as how many of them are real is often difficult to answer.

Non-backtracking operators first appeared in 1989 in the work of Hashimoto~\cite{hashimoto1989zeta}, in the context of zeta functions of graphs, and they have since proved useful in many areas, including algebraic graph theory \cite{terras2010zeta,bass1992ihara,cooper2009}, the study of Ramanujan graphings \cite{backhausz2015ramanujan}, epidemic spreading \cite{shrestha2015message,castellano2018relevance}, community detection \cite{krzakala2013spectral,bordenave2018nonbacktracking,bolla2025clustering}, and centrality measures \cite{martin2014localization,torres2021nonbacktracking}. We refer to \cite{torres2019non,coste2021eigenvalues,glover2021some,pastor2020localization,Arrigo1,Arrigo2,Gronford1,Mellor1,bolla2025non} for further results on non-backtracking operators and their applications. 

The idea behind these operators is to study walks on a graph that are never allowed to return along the edge they have just traversed. More precisely, let $G=(V,E)$ be a simple graph on $N$ nodes and $M$ edges, with minimum degree at least $2$. Each edge of $G$ admits two \emph{orientations}, and we write $[v,w]$ for the \emph{oriented edge} with \emph{input} $\inp([v,w])\coloneqq v$ and \emph{output} $\out([v,w])\coloneqq w$. Its inverse is $[v,w]^{-1}\coloneqq[w,v]$. This gives $2M$ oriented edges, which we label as $e_1,\ldots,e_{2M}$ in such a way that 
$$
e_{M+i}=e_i^{-1},\quad \text{for } i=1,\ldots,M.
$$

A \emph{non-backtracking random walk} on $G$ is a discrete-time Markov process on the oriented edges such that the probability of going from $e_i$ to $e_j$ is
\begin{equation*}
    \mathbb{P}(e_i\rightarrow e_j)=\begin{cases}\frac{1}{\deg(\out(e_i))-1} & \text{if }\out(e_i)=\inp(e_j)\text{ and }e_j\neq e_i^{-1},\\
    0 & \text{otherwise.}
    \end{cases}
\end{equation*}
In particular, the assumption that $G$ has minimum degree at least $2$ guarantees that the denominator is always positive. Moreover, one can also see a non-backtracking random walk as a process on $V$, but in this case the probability that the walker moves from a vertex $v$ to one of its neighbors depends on the vertex she came from. Hence, this process is not Markovian on $V$.

Non-backtracking random walks are encoded by the \emph{non-backtracking graph} or \emph{Hashimoto graph} of $G$, denoted $\mathcal{G}$, which is the directed graph whose vertices are the $2M$ oriented edges of $G$, and which has a directed edge from $e_i$ to $e_j$ if and only if $$\out(e_i)=\inp(e_j)  \text{ and } e_j\neq e_i^{-1},$$ i.e., if and only if a non-backtracking random walk can go from $e_i$ to $e_j$ in one step. The adjacency matrix of $\mathcal{G}$ is the $2M\times 2M$ matrix $\mathcal{A}$ with $(0,1)$--entries such that
\begin{equation*}
    \mathcal{A}_{ij}=1\iff \out(e_i)=\inp(e_j) \text{ and } e_j\neq e_i^{-1},
\end{equation*}
and its transpose $\mathcal{A}^\top$ is the \emph{non-backtracking matrix} of $G$ introduced by Hashimoto \cite{hashimoto1989zeta}. Moreover, a random walk on $\mathcal{G}$ is precisely a non-backtracking random walk on $G$. Since each vertex $e_i$ of $\mathcal{G}$ has (out-)degree $$\deg_{\mathcal{G}}(e_i)=\deg(\out(e_i))-1\geq 1,$$ the random walk Laplacian of $\mathcal{G}$ is well-defined, and it is called the \emph{non-backtracking Laplacian} of $G$. It can be written as
\begin{equation*}
    \mathcal{L}\coloneqq\id-\mathcal{D}^{-1}\mathcal{A}=\id-\mathcal{T},
\end{equation*}
where $\id$ denotes the identity matrix, $\mathcal{D}$ is the degree matrix of $\mathcal{G}$ and has diagonal entries $$\mathcal{D}_{ii}:=\deg_{\mathcal{G}}(e_i),$$ 
  and $\mathcal{T}\coloneqq\mathcal{D}^{-1}\mathcal{A}$ is the transition probability matrix of the non-backtracking random walk on $G$.

We introduced the non-backtracking Laplacian in \cite{JMT2023}, together with J.~Jost. Its spectral properties were further investigated in \cite{MZZ}, its connection to graph colorings was studied in \cite{maximal-colouring}, and it has been applied to vertex clustering in \cite{bolla2025clustering}. The computations in \cite[Table 2]{JMT2023} showed that, among graphs with at most ten nodes and minimum degree at least~$2$, substantially fewer graphs admit a cospectral mate for the non-backtracking Laplacian than for the other operators considered; that is, there are fewer pairs of non-isomorphic graphs with the same $\mathcal{L}$ spectrum. Moreover, in the same paper we showed that the presence of certain cycle configurations produces real eigenvalues of $\mathcal L$. 

For the non-backtracking matrix, the Ihara--Bass formula \cite{bass1992ihara,terras2010zeta,Kempton2016},
\begin{equation}\label{eq:Ihara}
    \det(\id-u\mathcal{A})=(1-u^2)^{M-N}\det\bigl(\id-uA+u^2(D-\id)\bigr),
\end{equation}
where $A$ and $D$ are the adjacency and degree matrices of $G$, gives at least $2M-2N$ real
eigenvalues for each of $\mathcal A$ and $\mathcal A^\top$,
counted with algebraic multiplicity. If $G$ is $d$--regular, then 
\[\mathcal L=\mathrm{Id}-\frac{1}{d-1}\mathcal A.\]
Therefore, in this case, $\mathcal L$ also has at least $2M-2N$
real eigenvalues, counted with algebraic multiplicity. For irregular graphs, however, degree normalization does not act by a scalar rescaling, and the Ihara--Bass formula does not directly give a corresponding lower bound for $\mathcal L$. It is therefore natural to ask whether the lower bound $2M-2N$
nevertheless holds for the non-backtracking Laplacian. 
Motivated by numerical experiments, we conjectured in 2022 that $\mathcal L$ always has at least $2M-2N$ real eigenvalues, counted with algebraic multiplicity. We formulated this conjecture during our course on \textit{non-backtracking operators of graphs} at the Max Planck Institute for Mathematics in the Sciences, and it was later stated as an open problem in \cite[Exercise~13]{Mulas2025LectureNotes}.

Here we prove this conjecture, and we refine it as follows.

\begin{theorem*}[Theorem~\ref{thm:refined} and Corollary~\ref{cor:sharp}]
Let $G$ be a simple graph on $N$ nodes and $M$ edges, with minimum degree at least $2$. Let $c$ be the number of connected components of $G$, and let $b$ be the number of bipartite ones. Then, the number of real eigenvalues of the non-backtracking Laplacian of $G$, counted with algebraic multiplicity, is at least
\begin{equation*}
    2M-2N+2c+2b.
\end{equation*}
Moreover, this bound is sharp.
\end{theorem*}

For connected regular graphs, we also determine exactly when the above bound is attained. Recall that a connected $d$--regular graph is \emph{Ramanujan} \cite{lubotzky1988ramanujan,murty2003ramanujan} if every adjacency eigenvalue $\mu\neq\pm d$ satisfies $$|\mu|\leq 2\sqrt{d-1}.$$  We call it \emph{strictly Ramanujan} if every adjacency eigenvalue $\mu\neq\pm d$ satisfies $$|\mu|<2\sqrt{d-1}.$$ 

We show that, for a connected regular graph $G$, the non-backtracking Laplacian has exactly $2M-2N+2+2b$ real eigenvalues if and only if $G$ is strictly Ramanujan. In particular, cycles, complete graphs and balanced complete bipartite graphs, which are strictly Ramanujan, attain the bound.

We also determine how these real eigenvalues are distributed around $1$. In Theorem~\ref{thm:split}, we show that at least $M-N+2c$ of them lie in $[0,1)$ and at least $M-N+2b$ lie in $(1,2]$. 

For completeness, we also prove in Theorem~\ref{thm:upper} that
$\mathcal L$ has at most $2M-2c-2b$ real eigenvalues, counted with algebraic multiplicity, and that this bound is attained by the cycle $C_4$.

\subsection*{Structure of the paper}

In Section~\ref{section:prelimnaries}, we fix the notation and establish some preliminary results, including Lemma~\ref{lemma:LA} from indefinite linear algebra. In Section~\ref{section:2M-2N}, we prove Theorem~\ref{thm:2022}, which settles our 2022 conjecture. In Section~\ref{section:refined}, we refine this bound by taking into account the connected and bipartite components of the graph, thus proving Theorem~\ref{thm:refined}. In Section~\ref{section:sharpness}, we show that the refined bound is sharp, and we prove that a connected regular graph attains it if and only if it is strictly Ramanujan. In Section~\ref{section:signs}, we show how the real eigenvalues split around $1$. Finally, in Section~\ref{section:upper}, we prove an upper bound on the number of real eigenvalues of the non-backtracking Laplacian.

\section{Notation and preliminaries}\label{section:prelimnaries}

Throughout the paper, we fix a simple graph $G=(V,E)$ on $N$ nodes and $M$ edges, with minimum degree at least $2$, and we use the same notation as in the introduction. In particular, $\mathcal{G}=(\mathcal{V},\mathcal{E})$ denotes the non-backtracking graph of $G$, whose vertices are the $2M$ oriented edges $e_1,\ldots,e_{2M}$ of $G$, labeled so that
$$
e_{M+i}=e_i^{-1},\quad \text{for } i=1,\ldots,M.
$$
Moreover, $\mathcal{A}$ and $\mathcal{D}$ denote the adjacency and degree matrices of $\mathcal{G}$, while $\mathcal{L}=\id-\mathcal{D}^{-1}\mathcal{A}$ and $\mathcal{T}=\mathcal{D}^{-1}\mathcal{A}$ denote the non-backtracking Laplacian and the transition matrix of the non-backtracking random walk on $G$.

We identify each vector $f\in\mathbb{C}^{2M}$ with the function
\begin{align*}
    f\colon\mathcal{V} &\rightarrow \mathbb{C}\\
    e_i&\mapsto f_i.
\end{align*}
 In these terms, the transition matrix acts as
\begin{equation}\label{eq:mathcalT}
    (\mathcal{T}f)([v,w])=\frac{1}{\deg w-1}\sum_{z\sim w,\,z\neq v}f([w,z]),
\end{equation}
for every $[v,w]\in \mathcal{V}$, where $z\sim w$ means that $z$ and $w$ are adjacent in $G$. 

Given a real symmetric matrix $K\in\mathbb{R}^{r\times r}$, the \emph{$K$-product} on $\mathbb{C}^r$ is defined as
\begin{equation*}
    (\mathbf{x},\mathbf{y})_K\coloneqq\langle\mathbf{x},K\mathbf{y}\rangle=\overline{\mathbf{x}}^\top K\mathbf{y},
\end{equation*}
where $\langle\cdot,\cdot\rangle$ is the standard inner product on $\mathbb{C}^r$. A matrix $R\in\mathbb{R}^{r\times r}$ is said to be \emph{self-adjoint with respect to the $K$-product}, or \emph{$K$-self-adjoint}, if $$(R\mathbf{x},\mathbf{y})_K=(\mathbf{x},R\mathbf{y})_K \text{ for all }\mathbf{x},\mathbf{y}\in\mathbb{C}^r,$$ or equivalently, if
\begin{equation*}
    R^\top K=KR.
\end{equation*}

As in \cite{JMT2023}, we let $P$ be the $2M\times 2M$ matrix
\begin{equation*}
    P\coloneqq\begin{pmatrix}0&\id\\ \id&0\end{pmatrix},
\end{equation*}
where each block has size $M\times M$. Since $e_{M+i}=e_i^{-1}$, we have that $$(Pf)(e)=f(e^{-1}) \text{ for each }e\in\mathcal{V}.$$ 
Moreover, by \cite[Theorem 3.8]{JMT2023}, $\mathcal{L}$ is $P$-self-adjoint, that is,
   \begin{equation}\label{eq:LP}
       \mathcal{L}^\top P=P\mathcal{L}.
   \end{equation}

We now introduce the matrix
\begin{equation*}
    H\coloneqq P(\id-\mathcal{L})=P\mathcal{T}.
\end{equation*} Its entries are
\begin{equation*}
    H_{ij}=\begin{cases}
        \frac{1}{\deg(\inp(e_i))-1} & \text{if } \inp(e_i)=\inp(e_j) \text{ and } e_i\neq e_j,\\
        0 & \text{otherwise.}
    \end{cases}
\end{equation*}
In particular, $H$ is real and symmetric. Moreover, $\mathcal{L}$ is $H$-self-adjoint. In fact, by \eqref{eq:LP} and since $\mathcal{L}$ commutes with $\id-\mathcal{L}$,
\begin{equation}\label{eq:LH}
    \mathcal{L}^\top H=\mathcal{L}^\top P(\id-\mathcal{L})=P\mathcal{L}(\id-\mathcal{L})=P(\id-\mathcal{L})\mathcal{L}=H\mathcal{L}.
\end{equation}

The matrix $H$ plays a central role in the proofs of our main results, for the following reason. By Lemma~\ref{lemma:LA} below, if a real matrix $R$ is self-adjoint with respect to the $K$-product of a real symmetric invertible matrix $K$, then $R$ has at most $2\min\bigl(n_+(K),n_-(K)\bigr)$ non-real eigenvalues, where $n_+(K)$ and $n_-(K)$ denote the numbers of positive and negative eigenvalues of $K$, counted with multiplicity. By \eqref{eq:LP}, $\mathcal{L}$ is $P$-self-adjoint. However, $P$ has eigenvalues $1$ and $-1$, each with multiplicity $M$, so in this case the lemma only says that $\mathcal{L}$ has at most $2M$ non-real eigenvalues, which is trivial since $\mathcal{L}$ is a $2M\times 2M$ matrix. The advantage of considering $H$ is that $\mathcal{L}$ is also $H$-self-adjoint, by \eqref{eq:LH}, and $H$ has only $N$ positive eigenvalues and $2M-N$ negative ones, by Lemma~\ref{lemma:H} below. Since every vertex of $G$ has degree at least $2$, we have that $$2M=\sum_{v\in V}\deg v\geq 2N,$$ and therefore 
\begin{equation}\label{eq:minN}
    \min\bigl(n_+(H),n_-(H)\bigr)=\min\bigl(N,2M-N\bigr)=N.
\end{equation}
Hence, by Lemma~\ref{lemma:LA}, $\mathcal{L}$ has at most $2N$ non-real eigenvalues, that is, at least $2M-2N$ real ones. This is the content of Theorem~\ref{thm:2022}. For the refinement in Theorem~\ref{thm:refined}, we then use the connected and bipartite components of $G$ to reduce the number of positive eigenvalues of $H$ that need to be taken into account.

We first prove the claim on the eigenvalues of $H$.

\begin{lemma}\label{lemma:H}
The matrix $H$ is invertible. Moreover,
\begin{equation*}
    n_+(H)=N\qquad\text{and}\qquad n_-(H)=2M-N.
\end{equation*}
\end{lemma}

\begin{proof}
By definition of $H$, we have that $H_{ij}=0$ whenever $\inp(e_i)\neq\inp(e_j)$. Hence, if we order the oriented edges according to their input, $H$ is block diagonal, with one block $H_v$ of size $\deg v\times\deg v$ for each $v\in V$. Moreover,
\begin{equation*}
    H_v=\frac{1}{\deg v-1}\bigl(J_{\deg v}-\id\bigr),
\end{equation*}
where $J_{\deg v}$ denotes the $\deg v\times\deg v$ matrix whose entries are all equal to $1$. Since $J_{\deg v}$ has eigenvalue $\deg v$ with multiplicity $1$, and eigenvalue $0$ with multiplicity $\deg v-1$, the block $H_v$ has eigenvalue $1$ with multiplicity $1$, and eigenvalue $-1/(\deg v-1)$ with multiplicity $\deg v-1$. The eigenvalues of $H$, with multiplicity, are those of the blocks $H_v$. Since none of them is zero, $H$ is invertible, and
\begin{equation*}
    n_+(H)=\sum_{v\in V}1=N,\qquad\text{while}\qquad n_-(H)=\sum_{v\in V}(\deg v-1)=2M-N.
\end{equation*}
\end{proof}

We now state the linear algebra fact used in the proofs of our main results. Given a real square matrix $R$, we let $n_{\mathbb{R}}(R)$ denote the number of real eigenvalues of $R$, counted with algebraic multiplicity.

\begin{lemma}\label{lemma:LA}
Let $K\in\mathbb{R}^{r\times r}$ be symmetric and invertible, and let $R\in\mathbb{R}^{r\times r}$ be $K$-self-adjoint. Then, $R$ has at most $2\min\bigl(n_+(K),n_-(K)\bigr)$ non-real eigenvalues, counted with algebraic multiplicity. Equivalently,
\begin{equation*}
    n_{\mathbb{R}}(R)\geq r-2\min\bigl(n_+(K),n_-(K)\bigr).
\end{equation*}
\end{lemma}

\begin{proof}
Since $R$ is $K$-self-adjoint, it follows from \cite[page 77]{gohberg2006indefinite} (see the inequality following equation (5.2.4) and the comment after it, where $A=R$, $H=K$ and $i_-(H)=n_-(K)$ in our notation) that $n_-(K)$ is at least half the number of non-real eigenvalues of $R$, counted with algebraic multiplicity. Since $R$ is also $(-K)$-self-adjoint, and since $n_-(-K)=n_+(K)$, the same holds for $n_+(K)$. Therefore, the number of non-real eigenvalues of $R$ is at most $2\min\bigl(n_+(K),n_-(K)\bigr)$.
\end{proof}

\section{The bound $2M-2N$}\label{section:2M-2N}

In this section, we prove our 2022 conjecture, which states that the non-backtracking Laplacian of a graph with minimum degree at least $2$ has at least $2M-2N$ real eigenvalues. As anticipated in Section~\ref{section:prelimnaries}, this follows by combining the $H$-self-adjointness of $\mathcal{L}$ with Lemma~\ref{lemma:H} and Lemma~\ref{lemma:LA}.

\begin{theorem}\label{thm:2022}
Let $G$ be a simple graph on $N$ nodes and $M$ edges, with minimum degree at least $2$. Then, the number of real eigenvalues of the non-backtracking Laplacian of $G$, counted with algebraic multiplicity, is at least
\begin{equation*}
    2M-2N.
\end{equation*}
\end{theorem}

\begin{proof} By \eqref{eq:LH}, $\mathcal{L}$ is $H$-self-adjoint, and by Lemma~\ref{lemma:H}, $H$ is a real symmetric invertible matrix with $n_+(H)=N$ and $n_-(H)=2M-N$. As observed in \eqref{eq:minN}, the minimum degree assumption implies that
\begin{equation*}
    \min\bigl(n_+(H),n_-(H)\bigr)=N.
\end{equation*}
Applying Lemma~\ref{lemma:LA} with $R=\mathcal{L}$ and $K=H$ yields
\begin{equation*}
    n_{\mathbb{R}}(\mathcal{L})\geq 2M-2N.\qedhere
\end{equation*}
\end{proof}

\section{The refined bound}\label{section:refined}

In this section, we refine Theorem~\ref{thm:2022} by taking into account the connected and bipartite components of $G$. These give explicit eigenvectors of $\mathcal{L}$ with real eigenvalues, which allow us to reduce the number of positive eigenvalues of $H$ that need to be taken into account in Lemma~\ref{lemma:LA}.

\begin{theorem}\label{thm:refined}
Let $G$ be a simple graph on $N$ nodes and $M$ edges, with minimum degree at least $2$. Let $c$ be the number of connected components of $G$, and let $b$ be the number of bipartite ones. Then, the number of real eigenvalues of the non-backtracking Laplacian of $G$, counted with algebraic multiplicity, is at least
\begin{equation*}
    2M-2N+2c+2b.
\end{equation*}
\end{theorem}
\begin{proof}
For each connected component $C$ of $G$, let $V(C)$ be the vertex set of $C$, let $E(C)$ be its edge set, and let $\mathbbm{1}_C\colon\mathcal{V}\to\{0,1\}$ be the indicator function of the oriented edges of $C$. If $C$ is bipartite, fix a function $s_C\colon V(C)\to\{1,-1\}$ that takes opposite values on adjacent vertices, and let
\begin{align*}
    \sigma_C\colon\mathcal{V}&\to\{0,1,-1\}\\
    [v,w]&\mapsto\begin{cases}
        s_C(v), & \text{if } [v,w] \text{ is an oriented edge of } C,\\
        0, & \text{otherwise.}
    \end{cases}
\end{align*}
By \eqref{eq:mathcalT} we have that, for every $[v,w]\in\mathcal{V}$,
\begin{align*}
    (\mathcal{T}\mathbbm{1}_C)([v,w])&=\frac{1}{\deg w-1}\sum_{z\sim w,\,z\neq v}\mathbbm{1}_C([w,z])\\
    &=\begin{cases}
        1 & \text{if } [v,w] \text{ is an oriented edge of } C,\\
        0 & \text{otherwise.}
    \end{cases}
\end{align*}
Hence, $\mathcal{T}\mathbbm{1}_C=\mathbbm{1}_C$. Similarly, for every $[v,w]\in\mathcal{V}$, by definition of $P$, and since $[v,w]$ and $[w,v]$ belong to the same connected component,
$$(P\mathbbm{1}_C)([v,w])=\mathbbm{1}_C([w,v])=\mathbbm{1}_C([v,w]).$$ 
 Hence, $P\mathbbm{1}_C=\mathbbm{1}_C$. Moreover, if $C$ is bipartite and $[v,w]$ is in $C$, then $s_C(w)=-s_C(v)$. Since each oriented edge $[w,z]$ with $z\sim w$ also belongs to $C$, this implies that, by \eqref{eq:mathcalT},
\begin{align*}
    (\mathcal{T}\sigma_C)([v,w])&=\frac{1}{\deg w-1}\sum_{z\sim w,\,z\neq v}\sigma_C([w,z])\\ &=\frac{1}{\deg w-1}\sum_{z\sim w,\,z\neq v}s_C(w)\\ &=-s_C(v),
\end{align*}
and by definition of $P$,
\begin{equation*}
    (P\sigma_C)([v,w])=\sigma_C([w,v])=s_C(w)=-s_C(v).
\end{equation*}
Together with the fact that both $\mathcal{T}\sigma_C$ and $P\sigma_C$ vanish outside $C$, this implies that
\begin{equation*}
    \mathcal{T}\sigma_C=-\sigma_C\qquad\text{and}\qquad P\sigma_C=-\sigma_C.
\end{equation*}
Since $H=P\mathcal{T}$ and $\mathcal{L}=\id-\mathcal{T}$, from the above equalities it follows that
\begin{equation*}
    H\mathbbm{1}_C=\mathbbm{1}_C,\qquad\mathcal{L}\mathbbm{1}_C=\mathbf{0},\qquad H\sigma_C=\sigma_C,\qquad\mathcal{L}\sigma_C=2\sigma_C.
\end{equation*}

Now, let $W\subseteq\mathbb{C}^{2M}$ be the span of the functions $\mathbbm{1}_C$, for all connected components $C$, and $\sigma_C$, for all bipartite components $C$. These $c+b$ functions are non-zero and mutually orthogonal with respect to the standard inner product. In fact, functions associated with different components have disjoint supports, and for a bipartite component $C$,
\begin{align*}
    \langle\mathbbm{1}_C,\sigma_C\rangle&=\sum_{e\in\mathcal{V}}\mathbbm{1}_C(e)\,\sigma_C(e)\\ &=\sum_{\{v,w\}\in E(C)}\bigl(\sigma_C([v,w])+\sigma_C([w,v])\bigr)\\ &=\sum_{\{v,w\}\in E(C)}\bigl(s_C(v)+s_C(w)\bigr)\\ &=0.
\end{align*}
Therefore, $\dim W=c+b$. Moreover, $H\mathbf{w}=\mathbf{w}$ for every $\mathbf{w}\in W$, and since $W$ is spanned by eigenvectors of $\mathcal{L}$, it is $\mathcal{L}$-invariant, that is, $\mathcal{L}\mathbf{w}\in W$ for every $\mathbf{w}\in W$.  

Now, for all $\mathbf{w}\in W$ and $\mathbf{x}\in\mathbb{C}^{2M}$, we have that
\begin{equation*}
    (\mathbf{w},\mathbf{x})_H=\overline{\mathbf{w}}^\top H\mathbf{x}=\bigl(H\overline{\mathbf{w}}\bigr)^\top \mathbf{x}=\overline{H\mathbf{w}}^\top \mathbf{x}=\overline{\mathbf{w}}^\top \mathbf{x}=\langle \mathbf{w},\mathbf{x}\rangle,
\end{equation*}
where the second equality holds because $H$ is symmetric, the third because $H$ is real, and the fourth because $H\mathbf{w}=\mathbf{w}$. 

Therefore, if we let $W^\perp$ be the orthogonal complement of $W$ with respect to the standard inner product, then $W^\perp$ is also $\mathcal{L}$-invariant. In fact, let $\mathbf{x}\in W^\perp$ and $\mathbf{w}\in W$. Then, since $\mathcal{L}\mathbf{w}\in W$, we have that
\begin{equation*}
    \langle\mathbf{w},\mathcal{L}\mathbf{x}\rangle=(\mathbf{w},\mathcal{L}\mathbf{x})_H=(\mathcal{L}\mathbf{w},\mathbf{x})_H=\langle\mathcal{L}\mathbf{w},\mathbf{x}\rangle=0,
\end{equation*}
where the first and third equalities follow from the above identity, the second one holds because $\mathcal{L}$ is $H$-self-adjoint, and the last one because $\mathbf{x}\in W^\perp$. Hence, $\mathcal{L}\mathbf{x}\in W^\perp$. 

Now, let $Q$ be a real orthogonal matrix whose first $c$ columns are the vectors $\mathbbm{1}_C/\|\mathbbm{1}_C\|$, whose next $b$ columns are the vectors $\sigma_C/\|\sigma_C\|$, and whose remaining $2M-c-b$ columns form an orthonormal basis of $W^\perp$, which can be chosen real since $W$ is spanned by real vectors. In other words, $Q$ is the change of basis to an orthonormal basis adapted to the decomposition $\mathbb{C}^{2M}=W\oplus W^\perp$. Since both summands are $\mathcal{L}$-invariant, $Q^\top\mathcal{L}Q$ is block diagonal, and since the first $c+b$ columns of $Q$ are eigenvectors of $\mathcal{L}$, with eigenvalue $0$ for the first $c$ ones and $2$ for the next $b$ ones,
\begin{equation*}
    Q^\top\mathcal{L}Q=\begin{pmatrix}\mathcal{L}_W&0\\0&\mathcal{L}_0\end{pmatrix},
\end{equation*}
where $\mathcal{L}_W$ is the diagonal matrix with $c$ diagonal entries equal to $0$ and $b$ diagonal entries equal to $2$, and $\mathcal{L}_0$ is a real matrix of size $(2M-c-b)\times(2M-c-b)$.

Similarly, let $\mathbf{q}_j$ denote the $j$-th column of $Q$. Then $(Q^\top HQ)_{ij}=(\mathbf{q}_i,\mathbf{q}_j)_H$, and by the identity above, $(\mathbf{q}_i,\mathbf{q}_j)_H=\langle\mathbf{q}_i,\mathbf{q}_j\rangle$, which is $1$ if $i=j$ and $0$ otherwise, whenever $\mathbf{q}_i\in W$ or $\mathbf{q}_j\in W$. Hence
\begin{equation*}
    Q^\top HQ=\begin{pmatrix}\id&0\\0&H_0\end{pmatrix},
\end{equation*}
where $H_0$ is a real symmetric matrix of size $(2M-c-b)\times(2M-c-b)$.

Since $H$ is invertible, so is $H_0$. Moreover, by Sylvester's law of inertia and Lemma~\ref{lemma:H},
\begin{equation*}
    n_+(H_0)=n_+(H)-(c+b)=N-c-b
\end{equation*}
and
\begin{equation*}
   n_-(H_0)=n_-(H)=2M-N.
\end{equation*}
Since $Q^\top=Q^{-1}$, the identity $\mathcal{L}^\top H=H\mathcal{L}$ implies that
\begin{equation*}
    (Q^\top\mathcal{L}Q)^\top(Q^\top HQ)=Q^\top\mathcal{L}^\top HQ=Q^\top H\mathcal{L}Q=(Q^\top HQ)(Q^\top\mathcal{L}Q),
\end{equation*}
and comparing the lower right blocks gives $\mathcal{L}_0^\top H_0=H_0\mathcal{L}_0$. That is, $\mathcal{L}_0$ is $H_0$-self-adjoint. By \eqref{eq:minN}, we have that $$N-c-b\leq N\leq 2M-N,$$ and therefore Lemma~\ref{lemma:LA} implies that $\mathcal{L}_0$ has at most $2(N-c-b)$ non-real eigenvalues, counted with algebraic multiplicity.

Finally, since $Q^\top\mathcal{L}Q$ is block diagonal,
\begin{equation*}
    \det(t\,\id-\mathcal{L})=\det(t\,\id-\mathcal{L}_W)\cdot\det(t\,\id-\mathcal{L}_0)=t^c(t-2)^b\det(t\,\id-\mathcal{L}_0).
\end{equation*}
Hence, every non-real eigenvalue of $\mathcal{L}$ is an eigenvalue of $\mathcal{L}_0$, with the same algebraic multiplicity, and therefore
\begin{equation*}
    n_{\mathbb{R}}(\mathcal{L})\geq 2M-2(N-c-b)=2M-2N+2c+2b.
\end{equation*}
\end{proof}

\section{Sharpness}\label{section:sharpness}

In this section, we show that the bound in Theorem~\ref{thm:refined} is sharp. For connected regular graphs, we first determine exactly when it is attained: by Theorem~\ref{thm:sharpness}, this happens if and only if the graph is strictly Ramanujan. As a consequence, in Corollary~\ref{cor:sharp} we obtain several explicit families of graphs attaining the bound.

\begin{theorem}\label{thm:sharpness}
    Let $G$ be a connected $d$-regular graph with $d\geq 2$, and let $\kappa$ be the number of adjacency eigenvalues $\mu\neq\pm d$ of $G$ with $|\mu|\geq 2\sqrt{d-1}$, counted with multiplicity. Then,
    \begin{equation*}
        n_{\mathbb{R}}(\mathcal{L})=2M-2N+2+2b+2\kappa.
    \end{equation*}
    In particular, equality in Theorem \ref{thm:refined} holds if and only if $G$ is strictly Ramanujan, that is, every adjacency eigenvalue $\mu\ne\pm d$ satisfies $|\mu|<2\sqrt{d-1}$.
\end{theorem}

\begin{proof}
Let $G$ be a connected $d$-regular graph on $N$ nodes and $M$ edges, with $d\geq 2$. Since $G$ is connected, we have that $c=1$, and therefore the quantity in Theorem~\ref{thm:refined} is
\begin{equation*}
    2M-2N+2+2b,
\end{equation*}
where $b=1$ if $G$ is bipartite, and $b=0$ otherwise.

Now, let $\mu_1,\ldots,\mu_N$ be the eigenvalues of the adjacency matrix $A$ of $G$, listed with multiplicity, and given $\mu\in\mathbb{R}$, let
\begin{equation*}
    q_\mu(t)\coloneqq t^2-\frac{\mu}{d-1}\,t+\frac{1}{d-1}.
\end{equation*}

Since $G$ is $d$-regular, we have that $D=d\,\id$ and
\begin{equation*}
    \mathcal{T}=\frac{1}{d-1}\mathcal{A}.
\end{equation*}
Hence, the Ihara--Bass formula \eqref{eq:Ihara} becomes
\begin{equation*}
    \det(\id-u\mathcal{A})=(1-u^2)^{M-N}\prod_{j=1}^N\bigl(1-\mu_ju+(d-1)u^2\bigr).
\end{equation*}
For $t\neq 0$, substituting $u=1/((d-1)t)$ and multiplying both sides by $t^{2M}$, we obtain
\begin{equation*}
    \det(t\,\id-\mathcal{T})=\left(t^2-\frac{1}{(d-1)^2}\right)^{M-N}\prod_{j=1}^{N}q_{\mu_j}(t).
\end{equation*}
Note that the exponent $M-N=N(d-2)/2$ is non-negative. We now count the real roots of each factor, with multiplicity.

The first factor has $2(M-N)$ real roots, namely $1/(d-1)$ and $-1/(d-1)$, each with multiplicity $M-N$.

Since $G$ is connected and $d$-regular, $d$ is a simple eigenvalue of $A$. Moreover, $-d$ is an eigenvalue of $A$ if and only if $G$ is bipartite, in which case it is also simple. The corresponding quadratic factors are
\begin{equation*}
    q_d(t)=(t-1)\left(t-\frac{1}{d-1}\right)\qquad\text{and}\qquad q_{-d}(t)=(t+1)\left(t+\frac{1}{d-1}\right),
\end{equation*}
which have two real roots each. Hence, the eigenvalues $\pm d$ contribute $2+2b$ real roots.

For every other eigenvalue $\mu\neq\pm d$ of $A$, the discriminant of $q_\mu$ is
\begin{equation*}
    \left(\frac{\mu}{d-1}\right)^2-\frac{4}{d-1}=\frac{\mu^2-4(d-1)}{(d-1)^2}.
\end{equation*}
Therefore, $q_\mu$ has two real roots, counted with multiplicity, if $|\mu|\geq 2\sqrt{d-1}$, and two non-real roots if $|\mu|<2\sqrt{d-1}$.

Since the multiplicities of the roots add up when polynomials are multiplied, we obtain, with $\kappa$ as in the statement,
\begin{equation*}
    n_{\mathbb{R}}(\mathcal{T})=2(M-N)+2+2b+2\kappa.
\end{equation*}
Since $\mathcal{L}=\id-\mathcal{T}$, this implies that
\begin{equation*}
    n_{\mathbb{R}}(\mathcal{L})=2M-2N+2+2b+2\kappa.
\end{equation*}
In particular, $\mathcal{L}$ has exactly $2M-2N+2+2b$ real eigenvalues if and only if $\kappa=0$, that is, if and only if $G$ is strictly Ramanujan.
\end{proof}

\begin{corollary}\label{cor:sharp}
    The bound in Theorem~\ref{thm:refined} is attained by the following graphs:
    \begin{itemize}
        \item the cycle graph $C_N$, for every $N\geq 3$;
        \item the complete graph $K_N$, for every $N\geq 3$;
        \item the complete bipartite graph $K_{r,r}$, for every $r\geq 2$.
    \end{itemize} In particular, the bound in Theorem~\ref{thm:refined} is sharp.
\end{corollary}

\begin{proof}
Since all these graphs are connected and regular, by Theorem~\ref{thm:sharpness} it suffices to show that they are strictly Ramanujan.
\begin{itemize}
    \item The cycle $C_N$ is $2$-regular, and its adjacency eigenvalues are $2\cos(2\pi j/N)$, for $j=0,\ldots,N-1$. All of them except $\pm 2$ satisfy $|\mu|<2=2\sqrt{2-1}$. 
\item The complete graph $K_N$ is $(N-1)$-regular, and its only adjacency eigenvalue other than $N-1$ is $-1$, which satisfies $1<2\sqrt{N-2}$ for $N\geq 3$. 
\item  The complete bipartite graph $K_{r,r}$ is $r$-regular, and its only adjacency eigenvalue other than $\pm r$ is $0$, which satisfies $0<2\sqrt{r-1}$ for $r\geq 2$. 
\end{itemize} This proves the claim.
\end{proof}

\section{Splitting the real eigenvalues around 1}\label{section:signs}

Theorem~\ref{thm:refined} counts the real eigenvalues of $\mathcal{L}$. In this section, we show how many of them must lie on each side of $1$. Equivalently, we count the positive and the negative real eigenvalues of $\mathcal{T}=\id-\mathcal{L}$. 

Recall from \cite[Theorem~4.1]{JMT2023} that every eigenvalue $\lambda$ of $\mathcal{L}$ satisfies $$|1-\lambda|\geq \frac{1}{\Delta-1},$$ where $\Delta$ is the largest degree of $G$. Hence, the real eigenvalues of $\mathcal{L}$ in $[0,1)$ actually lie in $[0,1-1/(\Delta-1)]$, and those in $(1,2]$ lie in $[1+1/(\Delta-1),2]$.

For a real symmetric matrix $K$, we define its \emph{signature} as $$\operatorname{sig}(K)\coloneqq n_+(K)-n_-(K).$$ For a real square matrix $R$, we let $n_>(R)$ and $n_<(R)$ denote the numbers of positive and negative real eigenvalues of $R$, counted with algebraic multiplicity.

The following lemma, which refines Lemma~\ref{lemma:LA}, is needed for the proof of Theorem~\ref{thm:split} below.

\begin{lemma}\label{lemma:signsplit}
Let $K\in\mathbb{R}^{r\times r}$ be symmetric and invertible, and let $R\in\mathbb{R}^{r\times r}$ be invertible and $K$-self-adjoint. Then,
\begin{equation*}
    n_>(R)\geq\frac{\bigl|\operatorname{sig}(K)+\operatorname{sig}(KR)\bigr|}{2}\qquad\text{and}\qquad n_<(R)\geq\frac{\bigl|\operatorname{sig}(K)-\operatorname{sig}(KR)\bigr|}{2}.
\end{equation*}
\end{lemma}

\begin{proof}
The matrix $KR$ is real, symmetric and invertible, since $$(KR)^\top=R^\top K=KR.$$ By the canonical form for $K$-self-adjoint matrices \cite[Theorem~5.1.1]{gohberg2006indefinite}, there is an invertible matrix $S$ such that $S^{-1}RS$ and $S^*KS$ are block diagonal with matching blocks, as follows. Each Jordan block $J_m(\lambda)$ of a real eigenvalue $\lambda$ is paired with $\varepsilon P_m$, where $P_m$ is the $m\times m$ matrix with ones on the anti-diagonal and zeros elsewhere, and $\varepsilon\in\{1,-1\}$. Each pair of Jordan blocks $J_m(\lambda)\oplus J_m(\overline{\lambda})$ of a non-real eigenvalue $\lambda$ is paired with $P_{2m}$. Since $$S^*KRS=(S^*KS)(S^{-1}RS),$$ the matrices $S^*KS$ and $S^*KRS$ are block diagonal, and by Sylvester's law of inertia their signatures are the sums of the signatures of the blocks.

The symmetric matrix $P_m$ has signature $1$ if $m$ is odd and $0$ if $m$ is even, and so does $P_mJ_m(\lambda)$ multiplied by $\operatorname{sgn}(\lambda)$, for real $\lambda\neq0$. For a pair of non-real eigenvalues, both blocks have signature $0$, since the span of the first $m$ basis vectors of the block is isotropic for both forms and has half the dimension. Therefore
\begin{equation*}
    \operatorname{sig}(K)=\sum\varepsilon,\qquad \operatorname{sig}(KR)=\sum\operatorname{sgn}(\lambda)\,\varepsilon,
\end{equation*}
where both sums run over the Jordan blocks of odd size belonging to real eigenvalues $\lambda$. Let $p$ and $q$ be the sums of $\varepsilon$ over such blocks with $\lambda>0$ and with $\lambda<0$, respectively. Then, $\operatorname{sig}(K)=p+q$ and $\operatorname{sig}(KR)=p-q$. Finally, $|p|$ is at most the number of blocks with $\lambda>0$, which is at most $n_>(R)$, and similarly $|q|\leq n_<(R)$.
\end{proof}

Using Lemma~\ref{lemma:signsplit}, we can now refine Theorem~\ref{thm:refined} as follows.

\begin{theorem}\label{thm:split}
Let $G$ be a simple graph on $N$ nodes and $M$ edges, with minimum degree at least $2$. Let $c$ be the number of connected components of $G$, and let $b$ be the number of bipartite ones. Then, $\mathcal{T}$ has at least $M-N+2c$ positive and at least $M-N+2b$ negative real eigenvalues, counted with algebraic multiplicity. Equivalently, $\mathcal{L}$ has at least $M-N+2c$ real eigenvalues in $[0,1)$ and at least $M-N+2b$ real eigenvalues in $(1,2]$.
\end{theorem}

\begin{proof}
We use the notation of the proof of Theorem~\ref{thm:refined}. Let $\mathcal{T}_0\coloneqq \id-\mathcal{L}_0$, so that $Q^\top\mathcal{T}Q$ is block diagonal with blocks $\id-\mathcal{L}_W$ and $\mathcal{T}_0$, and $Q^\top HQ=Q^\top P\mathcal{T}Q$ has blocks $\id$ and $H_0$.

Since $P\mathbbm{1}_C=\mathbbm{1}_C$ and $P\sigma_C=-\sigma_C$, the matrix $P$ maps $W$ to itself, and since $P$ is symmetric, it also maps $W^\perp$ to itself. Hence $Q^\top PQ$ is block diagonal, with a diagonal first block $P_W$ having $c$ entries equal to $1$ and $b$ entries equal to $-1$, and a second block $P_0$. Comparing the second blocks in $$Q^\top HQ=(Q^\top PQ)(Q^\top\mathcal{T}Q)$$ and in \eqref{eq:LP} gives
\begin{equation*}
    H_0=P_0\mathcal{T}_0\qquad\text{and}\qquad \mathcal{T}_0^\top P_0=P_0\mathcal{T}_0.
\end{equation*}
The matrix $P$ has eigenvalues $1$ and $-1$, each with multiplicity $M$. By Sylvester's law of inertia, $n_+(P_0)=M-c$ and $n_-(P_0)=M-b$, and from the proof of Theorem~\ref{thm:refined}, $n_+(H_0)=N-c-b$ and $n_-(H_0)=2M-N$. Hence,
\begin{equation*}
    \operatorname{sig}(P_0)=b-c\qquad\text{and}\qquad\operatorname{sig}(H_0)=2N-2M-c-b.
\end{equation*}
Since $H_0$ is invertible, so is $\mathcal{T}_0$. Lemma~\ref{lemma:signsplit} with $K=P_0$ and $R=\mathcal{T}_0$ gives
\begin{equation*}
    n_>(\mathcal{T}_0)\geq\frac{|(b-c)+(2N-2M-c-b)|}{2}=M-N+c
\end{equation*}
and
\begin{equation*}
  n_<(\mathcal{T}_0)\geq\frac{|(b-c)-(2N-2M-c-b)|}{2}=M-N+b,
\end{equation*}
where we used $M\geq N$. Finally, the first block $\id-\mathcal{L}_W$ contributes the eigenvalue $1$ with multiplicity $c$ and the eigenvalue $-1$ with multiplicity $b$. The statement for $\mathcal{L}$ follows from $\mathcal{L}=\id-\mathcal{T}$ and the fact that all eigenvalues of $\mathcal{T}$ have modulus at most $1$, since $\mathcal{T}$ is row-stochastic.
\end{proof}

Adding the two bounds of Theorem~\ref{thm:split} gives Theorem~\ref{thm:refined}. Both bounds are attained by the strictly Ramanujan graphs of Theorem~\ref{thm:sharpness}: in the proof of that theorem, the factor $$\bigl(t^2-1/(d-1)^2\bigr)^{M-N}$$ contributes $M-N$ positive and $M-N$ negative roots, the eigenvalue $d$ of $A$ contributes the positive roots $1$ and $1/(d-1)$, and the eigenvalue $-d$, present exactly when $G$ is bipartite, contributes the negative roots $-1$ and $-1/(d-1)$.

\begin{remark}
The results of this section do not depend on the particular weights $1/(\deg v-1)$, and they hold for any weighted matrix $\mathcal{T}_f\coloneqq C_f\mathcal{A}$, where $f$ is a positive function on the vertices and $C_f$ is the diagonal matrix with diagonal entries $(C_f)_{ee}\coloneqq f(\out(e))$. In fact, as in Lemma~\ref{lemma:H}, $P\mathcal{T}_f$ is symmetric with $N$ positive and $2M-N$ negative eigenvalues, therefore Lemma~\ref{lemma:signsplit} shows that $\mathcal{T}_f$ has at least $M-N$ positive and at least $M-N$ negative real eigenvalues. The refinement by $2c$ and $2b$ can also be obtained, using the Perron--Frobenius Theorem in place of the explicit eigenvectors $\mathbbm{1}_C$ and $\sigma_C$. In particular, $\mathcal{T}$ corresponds to $f(v)=1/(\deg v-1)$, while $f\equiv 1$ gives the transpose of the non-backtracking matrix. Note, however, that the location of the real eigenvalues of $\id-\mathcal{T}$ in $[0,1)$ and $(1,2]$ uses the fact that $\mathcal{T}$ is row-stochastic, and it does not extend to arbitrary weights.
\end{remark}

\section{An upper bound}\label{section:upper}

We now complement Theorem~\ref{thm:refined} with an upper bound on the number of real eigenvalues of $\mathcal{L}$.

\begin{theorem}\label{thm:upper}
 Let $G$ be a simple graph on $N$ nodes and $M$ edges, with minimum degree at least $2$. Let $c$ be the number of connected components of $G$, and let $b$ be the number of bipartite ones. Then, the number of real eigenvalues of the non-backtracking Laplacian of $G$, counted with algebraic multiplicity, is at most
\begin{equation*}
    2M-2c-2b.
\end{equation*} Moreover, this bound is attained by the cycle $C_4$.
\end{theorem}

\begin{proof}

Since $\mathcal{L}=\id-\mathcal{T}$, it suffices to show that $\mathcal{T}$ has at least $2c+2b$ non-real eigenvalues, counted with algebraic multiplicity.

If we order the oriented edges of $G$ according to the connected components of $G$, then $\mathcal{T}$ is block diagonal, with one block $\mathcal{T}_C$ for each connected component $C$, and $\mathcal{T}_C$ is the transition matrix of the non-backtracking random walk on $C$. Hence, it suffices to show that each $\mathcal{T}_C$ has at least two non-real eigenvalues, and at least four if $C$ is bipartite. To this end, we use the following three facts.

\begin{enumerate}[label=(\roman*)]
    \item All eigenvalues of $\mathcal{T}_C$ are non-zero. In fact, $\mathcal{T}=PH$ is invertible, since $P^2=\id$ and $H$ is invertible by Lemma~\ref{lemma:H}, and the eigenvalues of $\mathcal{T}_C$ are eigenvalues of $\mathcal{T}$.
    \item Since $\mathcal{T}_C$ is row-stochastic, all eigenvalues $t$ of $\mathcal{T}_C$ satisfy $|t|\leq 1$. Moreover, for the same reason, $1$ is an eigenvalue of $\mathcal{T}_C$, with the constant functions as eigenfunctions.
    \item We have that $\operatorname{tr}(\mathcal{T}_C^2)=0$. In fact,
$$
\operatorname{tr}(\mathcal{T}_C^2)=\sum_{e}(\mathcal{T}_C^2)_{ee}=\sum_{e}\sum_{e'}(\mathcal{T}_C)_{ee'}(\mathcal{T}_C)_{e'e},
$$ where $e$ and $e'$ range over the oriented edges of $C$. If $(\mathcal{T}_C)_{ee'}(\mathcal{T}_C)_{e'e}\neq 0$, then both factors are non-zero, that is, a non-backtracking random walk can go from $e$ to $e'$ and from $e'$ to $e$. The first step requires $e'$ to start at the output of $e$, and the second one requires $e'$ to end at the input of $e$, therefore $e'=e^{-1}$. But the step from $e$ to $e^{-1}$ is excluded by the definition of a non-backtracking random walk. Hence, all terms in the above sum vanish.
\end{enumerate}
Now, let $t_1,\ldots,t_m$ be the eigenvalues of $\mathcal{T}_C$, listed with algebraic multiplicity. Then, by (iii),
\begin{equation}\label{eq:trace}
    \sum_{k=1}^m t_k^2=\operatorname{tr}(\mathcal{T}_C^2)=0.
\end{equation}

By (i), every real eigenvalue $t_k$ of $\mathcal{T}_C$ satisfies $t_k^2>0$. Hence, by \eqref{eq:trace}, $\mathcal{T}_C$ has at least one non-real eigenvalue. Since $\mathcal{T}_C$ is real, its non-real eigenvalues appear in complex conjugate pairs, and therefore $\mathcal{T}_C$ has at least two non-real eigenvalues. 

Now, assume that $C$ is bipartite. Let $V(C)$ denote the vertex set of $C$, and let $M_C$ denote the number of edges of $C$. We want to show that, in this case, $\mathcal{T}_C$ cannot have exactly one pair of complex conjugate non-real eigenvalues. Since $C$ is bipartite, there exists a function $s\colon V(C)\to\{1,-1\}$ that takes opposite values on adjacent vertices. Hence, we can define the $2M_C\times 2M_C$ diagonal matrix $Z_C$ with diagonal entries
\begin{equation*}
    (Z_C)_{ee}\coloneqq s(\inp(e)),
\end{equation*}
for each oriented edge $e$ of $C$. 

If $(\mathcal{T}_C)_{ee'}\neq 0$, then a non-backtracking random walk can go from $e$ to $e'$, implying that $\inp(e')=\out(e)$, and therefore $s(\inp(e'))=-s(\inp(e))$. Therefore, since $Z_C$ is diagonal, for all oriented edges $e,e'$ of $C$ we have that
\begin{align*}
    (Z_C\mathcal{T}_CZ_C)_{ee'}&=\sum_{f,g}(Z_C)_{ef}\,(\mathcal{T}_C)_{fg}\,(Z_C)_{ge'}\\
    &=(Z_C)_{ee}\,(\mathcal{T}_C)_{ee'}\,(Z_C)_{e'e'}\\
    &=s(\inp(e))\,(\mathcal{T}_C)_{ee'}\,s(\inp(e'))\\
    &=-(\mathcal{T}_C)_{ee'},
\end{align*}
where the last equality follows from the previous observation if $(\mathcal{T}_C)_{ee'}\neq 0$, and is trivial if $(\mathcal{T}_C)_{ee'}=0$. Hence, $$Z_C\mathcal{T}_CZ_C=-\mathcal{T}_C.$$
Since $Z_C^2=\id$, it follows that $\mathcal{T}_C$ and $-\mathcal{T}_C$ are similar, and therefore the eigenvalues of $\mathcal{T}_C$, with multiplicity, are invariant under $t\mapsto -t$. In particular, since $1$ is an eigenvalue of $\mathcal{T}_C$ by (ii), so is $-1$. 

Now assume, by contradiction, that $\mathcal{T}_C$ has exactly two non-real eigenvalues, $t$ and $\overline{t}$. By the above invariance, $-t$ is also a non-real eigenvalue of $\mathcal{T}_C$, therefore $-t\in\{t,\overline{t}\}$. Since $t\neq 0$ by (i), we have that $-t\neq t$, and therefore $-t=\overline{t}$. This means that $t$ is purely imaginary, that is, $t=i\gamma$ for some $\gamma\in\mathbb{R}$ with $\gamma\neq 0$. Moreover, $|\gamma|=|t|\leq 1$ by (ii).

We now separate, in \eqref{eq:trace}, one copy of the eigenvalue $1$, one copy of the eigenvalue $-1$, and the two eigenvalues $i\gamma$ and $-i\gamma$. Denoting by $I$ the set of indices of the remaining eigenvalues, which are all real, we obtain
\begin{equation*}
    0=\sum_{k=1}^m t_k^2=1^2+(-1)^2+(i\gamma)^2+(-i\gamma)^2+\sum_{k\in I}t_k^2=2-2\gamma^2+\sum_{k\in I}t_k^2.
\end{equation*} 
Since $|\gamma|\leq 1$, we have that $2-2\gamma^2\geq 0$, and by (i), $t_k^2>0$ for every $k\in I$. Hence, both $2-2\gamma^2$ and the sum over $I$ must vanish, that is, $\gamma^2=1$ and $I=\emptyset$. Therefore, the only eigenvalues of $\mathcal{T}_C$ are $1$, $-1$, $i$ and $-i$, each with multiplicity one. Since $\mathcal{T}_C$ is a $2M_C\times 2M_C$ matrix, this implies that $M_C=2$. However, since every vertex of $C$ has degree at least $2$, we have $$2M_C=\sum_{v\in V(C)}\deg v\geq 2|V(C)|,$$ implying that $M_C\geq|V(C)|\geq 3$. This is a contradiction, and therefore $\mathcal{T}_C$ has at least four non-real eigenvalues.

Finally, for the cycle $C_4$, we have that $M=N=4$ and $c=b=1$, therefore
\begin{equation*}
    2M-2N+2c+2b=4=2M-2c-2b.
\end{equation*}
Together with Theorem~\ref{thm:refined}, this shows that $\mathcal{L}$ has exactly $4$ real eigenvalues, and that the bound is attained.
\end{proof}

\section*{Declaration of AI assistance}
The proofs presented in this paper were developed with the assistance of ChatGPT-6.1 Sol (OpenAI) and Claude Opus 5.5 (Anthropic). All proofs were checked in detail by the authors, who take full responsibility for the content.

\section*{Funding}
Raffaella Mulas is supported by the Dutch Research Council (NWO) through the grant VI.Veni.232.002 (Grant ID: \url{https://doi.org/10.61686/NLUFV33835}).

\bibliographystyle{plain} 

\bibliography{Bibliography}

@article{Kempton2016,
  title={Non-Backtracking Random Walks and a Weighted {I}hara’s Theorem},
  author={Kempton, M.},
  journal={Open Journal of Discrete Mathematics},
  volume={6},
  number={4},
  pages={207--226},
  year={2016},
  publisher={Scientific Research Publishing}
}

@book{gohberg2006indefinite,
  title={Indefinite linear algebra and applications},
  author={Gohberg, I. and Lancaster, P. and Rodman, L.},
  year={2006},
  publisher={Springer Science \& Business Media}
}

@book{chung,
  title={Spectral graph theory},
  author={Chung, F.},  year={1997},
  publisher={American Mathematical Soc.}
}

@book{brouwer,
  title={Spectra of graphs},
  author={Brouwer, A.E. and Haemers, W.H.},
  year={2011},
  publisher={Springer Science \& Business Media}
}

@article{torres2021nonbacktracking,
  title={Nonbacktracking eigenvalues under node removal: X-centrality and targeted immunization},
  author={Torres, L. and Chan, K.S. and Tong, H. and Eliassi-Rad, T.},
  journal={SIAM Journal on Mathematics of Data Science},
  volume={3},
  number={2},
  pages={656--675},
  year={2021},
  publisher={SIAM}
}

@article{bordenave2018nonbacktracking,
  title={{Nonbacktracking spectrum of random graphs: Community detection and nonregular Ramanujan graphs}},
  author={Bordenave, C. and Lelarge, M. and Massouli{\'e}, L.},
  journal={Annals of probability: An official journal of the Institute of Mathematical Statistics},
  volume={46},
  number={1},
  pages={1--71},
  year={2018},
  publisher={Institute of Mathematical Statistics}
}

@book{terras2010zeta,
  title={Zeta functions of graphs: {A} stroll through the garden},
  author={Terras, A.},
  series={{C}ambridge studies in advanced mathematics},
  volume={128},
  year={2010},
  publisher={Cambridge University Press}
}

@incollection{hashimoto1989zeta,
  title={Zeta functions of finite graphs and representations of $p$-adic groups},
  author={Hashimoto, K.},
  booktitle={Automorphic Forms and Geometry of Arithmetic Varieties},
  series={{A}dvanced {S}tudies in {P}ure {M}athematics},
  volume=15,
  pages={211--280},
  year={1989},
  publisher={Elsevier}
}

@article{bass1992ihara,
  title={{The Ihara-Selberg zeta function of a tree lattice}},
  author={Bass, H.},
  journal={International Journal of Mathematics},
  volume={3},
  number={06},
  pages={717--797},
  year={1992},
  publisher={World Scientific}
}

@article{krzakala2013spectral,
  title={Spectral redemption in clustering sparse networks},
  author={Krzakala, F. and Moore, C. and Mossel, E. and Neeman, J. and Sly, A. and Zdeborov{\'a}, L. and Zhang, P.},
  journal={Proceedings of the National Academy of Sciences},
  volume={110},
  number={52},
  pages={20935--20940},
  year={2013},
  publisher={National Acad Sciences}
}

@article{martin2014localization,
  title={Localization and centrality in networks},
  author={Martin, T. and Zhang, X. and Newman, M.E.J.},
  journal={Physical Review E},
  volume={90},
  number={5},
  pages={052808},
  year={2014},
  publisher={APS}
}

@article{shrestha2015message,
  title={Message-passing approach for recurrent-state epidemic models on networks},
  author={Shrestha, M. and Scarpino, S.V. and Moore, C.},
  journal={Physical Review E},
  volume={92},
  number={2},
  pages={022821},
  year={2015},
  publisher={APS}
}

@article{castellano2018relevance,
  title={Relevance of backtracking paths in recurrent-state epidemic spreading on networks},
  author={Castellano, C. and Pastor-Satorras, R.},
  journal={Physical Review E},
  volume={98},
  number={5},
  pages={052313},
  year={2018},
  publisher={APS}
}

@misc{Mulas2025LectureNotes,
  author       = {Mulas, R.},
  title        = {Spectral Theory of Graphs and Hypergraphs},
  howpublished = {Lecture notes for the BeyondTheEdge MSCA Doctoral Network,
                  Marseille, France},
  year         = {2025},
  note         = {Available at \url{https://cdn.prod.website-files.com/673b1c6a31679f06b3cef571/68d3daeb1b89a52be5dee00e_Marseille-LectureNotes-RM.pdf} (Accessed October 2026)}
}

@article{JMT2023,
  title={{Spectral theory of the non-backtracking Laplacian for graphs}},
  author={Jost, J. and Mulas, R. and Torres, L.},
  journal={Discrete Mathematics},
  year={2023},
  volume={346},
  number={10},
  pages={113536}
}

@book{BookJMZ,
place={Cambridge}, series={Cambridge Studies in Advanced Mathematics}, 
title={Spectra of Discrete Structures}, publisher={Cambridge University Press}, author={Jost, J. and Mulas, R. and Zhang, D.}, year={2026}, 
number ={221},
collection={Cambridge Studies in Advanced Mathematics}}

@article{backhausz2015ramanujan,
  title={Ramanujan graphings and correlation decay in local algorithms},
  author={Backhausz, {\'A}. and Szegedy, B. and Vir{\'a}g, B.},
  journal={Random Structures \& Algorithms},
  volume={47},
  number={3},
  pages={424--435},
  year={2015},
  publisher={Wiley Online Library}
}

@article{cooper2009,
  title={{Properties determined by the Ihara zeta function of a graph}},
  author={Cooper, Y.},
  journal={The Electronic Journal of Combinatorics},
  volume={16},
  number={1},
  pages={R84},
  year={2009}
}

@article{torres2019non,
  title={Non-backtracking cycles: length spectrum theory and graph mining applications},
  author={Torres, L. and Su{\'a}rez-Serrato, P. and Eliassi-Rad, T.},
  journal={Applied Network Science},
  volume={4},
  number={1},
  pages={1--35},
  year={2019},
  publisher={Springer}
}

@article{coste2021eigenvalues,
  title={Eigenvalues of the non-backtracking operator detached from the bulk},
  author={Coste, S. and Zhu, Y.},
  journal={Random Matrices: Theory and Applications},
  volume={10},
  number={03},
  pages={2150028},
  year={2021},
  publisher={World Scientific}
}

@article{glover2021some,
  title={Some spectral properties of the non-backtracking matrix of a graph},
  author={Glover, C. and Kempton, M.},
  journal={Linear Algebra and its Applications},
  volume={618},
  pages={37--57},
  year={2021},
  publisher={Elsevier}
}

@article{pastor2020localization,
  title={The localization of non-backtracking centrality in networks and its physical consequences},
  author={Pastor-Satorras, R. and Castellano, C.},
  journal={Scientific reports},
  volume={10},
  number={1},
  pages={1--12},
  year={2020},
  publisher={Nature Publishing Group}
}

@article{Arrigo1,
  title={Beyond non-backtracking: non-cycling network centrality measures},
  author={Arrigo, F. and Higham, D.J. and Noferini, V.},
  journal={Proceedings. Mathematical, Physical, and Engineering Sciences},
  year={2020},
  volume={476}
}

@article{Arrigo2,
  title={Weighted enumeration of nonbacktracking walks on weighted graphs},
   author={Arrigo, F. and Higham, D.J. and Noferini, V. and Wood, R.},
  journal={SIAM Journal on Matrix Analysis and Applications},
  volume={45},
  number={1},
  pages={397--418},
  year={2024},
  publisher={SIAM}
}

@article{Gronford1,
author = {Grindrod, P. and Higham, D.J. and Noferini, V.},
title = {{The Deformed Graph Laplacian and Its Applications to Network Centrality Analysis}},
journal = {SIAM Journal on Matrix Analysis and Applications},
volume = {39},
number = {1},
pages = {310-341},
year = {2018},
doi = {10.1137/17M1112297},
URL = {    https://doi.org/10.1137/17M1112297}
    }

@article{Mellor1,
  title = {Graph comparison via the nonbacktracking spectrum},
  author = {Mellor, A. and Grusovin, A.},
  journal = {Phys. Rev. E},
  volume = {99},
  issue = {5},
  pages = {052309},
  numpages = {7},
  year = {2019},
  month = {May},
  publisher = {American Physical Society},
  doi = {10.1103/PhysRevE.99.052309},
  url = {https://link.aps.org/doi/10.1103/PhysRevE.99.052309}
}

@article{bolla2025non,
  title={The non-backtracking random walk and its usage for vertex clustering},
  author={Bolla, M.},
  journal={arXiv preprint arXiv:2512.24434},
  year={2025}
}

@inproceedings{bolla2025clustering,
  title={Clustering the nodes of sparse edge-weighted graphs via non-backtracking spectra},
  author={Bolla, M. and Reittu, H. and Abdelkhalek, F.},
  booktitle={Future of Information and Communication Conference},
  pages={81--99},
  year={2025},
  organization={Springer}
}

@article{MZZ,
  title={{There is no going back: Properties of the non-backtracking Laplacian}},
  author={Mulas, R. and Zhang, D. and Zucal, G.},
  journal={Linear Algebra Appl.},
  VOLUME = {680},
      YEAR = {2024},
     PAGES = {341--370}
}

@article{maximal-colouring,
  title={{Maximal colourings for graphs}},
  author={Mulas, R.},
  journal={Graphs and Combinatorics},
  volume={40},
  number={93},
doi = {10.1007/s00373-024-02823-3},
  year={2024}
}

@article{lubotzky1988ramanujan,
  title={Ramanujan graphs},
  author={Lubotzky, A. and Phillips, R. and Sarnak, P.},
  journal={Combinatorica},
  volume={8},
  number={3},
  pages={261--277},
  year={1988},
  publisher={Springer-Verlag Berlin/Heidelberg}
}

@article{murty2003ramanujan,
  title={Ramanujan graphs},
  author={Murty, M.R.},
  journal={Journal of the Ramanujan Mathematical Society},
  volume={18},
  number={1},
  pages={33--52},
  year={2003},
  publisher={THE RAMANUJAN MATHEMATICAL SOCIETY}
}

\end{document}